\documentclass[11pt,twoside]{amsart}
\usepackage{amscd}

\newtheorem{theorem}{Theorem}[section]

\newtheorem{lemma}[theorem]{Lemma}

\newtheorem{proposition}[theorem]{Proposition}
\theoremstyle{definition}

\theoremstyle{remark}
\newtheorem{remark}[theorem]{Remark}

\numberwithin{equation}{section}

\numberwithin{equation}{subsection}
\newcommand{\be}
  {\protect\setcounter{equation}{\value{subsubsection}}}  
  \newcommand{\ee}%
   {\protect\setcounter{subsubsection}{\value{equation}}}

  {\protect\setcounter{subsubsection}{\value{equation}}}

\def \bC{\mathbb C}
\def \codim{\operatorname{codim}}

\def \diag{\operatorname{diag}}

\def \bG{\mathbb G}

\def \GL{\operatorname{GL}}

\def \bH{\mathbb H}
\def \cH{\mathcal H}

\def \ker{\operatorname{ker}}

\def \cO{\mathcal O}

\def \PGL{\operatorname{PGL}}

\def \bP{\mathbb P}

\def \bQ{\mathbb  Q}

\def \rk{\operatorname{rk}}
\def \res{respectively}

\def \bR{\mathbb R}

\def \SL{\operatorname{SL}}
\def \cS{\mathcal S}

\def \bZ{\mathbb Z}

\begin{document}

\title{Intersection cohomology of reductive varieties} 
\author{Michel Brion and Roy Joshua}
\address{Institut Fourier, B.P. 74, 38402 Saint Martin d'H\`eres Cedex,
France.}
\email{Michel.Brion@ujf-grenoble.fr}

\address{Department of Mathematics, Ohio State University, Columbus,
Ohio, 43210, USA.}
\curraddr{School of Mathematics, Institute for Advanced Study, Princeton, New Jersey, 08540, USA}
\email{joshua@ias.edu, joshua@math.ohio-state.edu}

\thanks{The second author is supported by the IHES, MPI (Bonn) and a
grant from the National Security Agency}



 \catcode`\@=11
 \def\leftheadline{\rlap{\ }\hfill\iftrue\topmark\fi\hfill}
 \def\rightheadline{\hfill\expandafter\iffalse\botmark\fi\hfill\llap{\ }}
 \def\output@{\shipout\vbox{%
  \iffirstpage@ \global\firstpage@false
  \pagebody \logo@ \makefootline%
 \else \ifrunheads@ \makeheadline \pagebody \makefootline
        \else \pagebo \makefootline \fi
 \fi}%
  \advancepageno \ifnum\outputpenalty>-\@M\else\dosupereject\fi}
 \catcode`\@=\active
 
\begin{abstract}
We extend the methods developed in our earlier work to algorithmically
compute the intersection cohomology Betti numbers of 
reductive varieties. These form a class of highly symmetric varieties
that includes equivariant compactifications of reductive groups. 
Thereby, we extend a well-known algorithm for toric varieties. 
\end{abstract}

\maketitle

\section{Introduction}

In this paper we extend the methods developed in our previous work
\cite{BJ} to {\it algorithmically} compute the local and global
intersection cohomology Betti numbers of a large class of varieties
with group action, that includes toric varieties.

\medskip

Let $G$ be a complex connected reductive algebraic group and let $B$
be a Borel subgroup. A normal complex algebraic variety $X$, equipped
with an action of $G$, is {\it spherical} if it contains a dense orbit
of $B$. We say that $X$ is {\it scs} (simply-connected spherical) if,
in addition, the $B$-isotropy group of this dense orbit is
connected. Equivalently, any $B$-equivariant finite surjective
morphism from a $B$-variety to $X$ is an isomorphism.

\medskip

For example, the group $G$ is a scs variety with respect to the action
of $G\times G$ by left and right multiplication ; thus, all normal
$G\times G$-equivariant embeddings of $G$ are scs as well. In
particular, all toric varieties are scs. Other examples include  
the space of all skew bilinear forms in $n$ variables for the
natural action of $\GL_n$ or $\SL_n$, and its subvarieties of
forms of rank at most $r$. But the space of symmetric
bilinear forms in $n$ variables is spherical, not scs.

\medskip

Spherical $G$-varieties enjoy the following properties:

\begin{itemize}
\item they contain only finitely many $B$-orbits, and hence only
finitely many $G$-orbits,
\item each $G$-orbit admits a {\it slice} (see ~\ref{slice} below) which
is an affine spherical variety under a connected reductive subgroup of
$G$, and 
\item the associated {\it link} (see ~\ref{att.fixed.pt} below) is
a projective spherical variety, of strictly smaller dimension.
\end{itemize}

\medskip

Important additional features of scs varieties are

\begin{itemize}
\item the connectedness of $B$- and $G$-isotropy groups of all points,
and 
\item the fact that all slices and links are scs as well. 
\end{itemize}

\medskip

This makes scs varieties particularly suited for applying the methods
from \cite{BJ}. They yield two recursive relations, expressing the
global intersection cohomology Betti numbers of projective scs
varieties in terms of the corresponding local numbers, and the latter
in terms of the global numbers of the links. This is the content of
our main result.

\begin{theorem}\label{main} 
Let $X$ be a projective scs $G$-variety; let $IP_X(t)$
($IP_{X,x}(t)$) be the Poincar\'e polynomial for global intersection
cohomology (for the stalks of the intersection cohomology
sheaves at $x\in X$, \res). Then 

\be\begin{equation}\label{A.1}
IP_X(t) = \sum_x  (t^2 - 1)^{r-r_x} 
\frac{P_{G/T}(t)}{P_{G_x/T_x}(t)} IP_{X,x}(t)
\end{equation}\ee
(sum over representatives of $G$-orbits in $X$), and 
\be\begin{equation}\label{A.2}
IP_{X,x}(t) = IP_{\cS_x,x}(t) = 
\tau_{\le d_x-1}((1-t^2)IP_{\bP(S_x)}(t)). 
\end{equation} \ee
Here $T$ is a maximal torus of $G$, of dimension $r$ ; $T_x$ is
a maximal torus of $G_x$, of dimension $r_x$; $\cS_x$ is a slice
at $x$ to $Gx$, of dimension $d_x$;  
$\bP(\cS_x) = (S_x-x)/\bG_m$ is the corresponding link for an
attractive action of the multiplicative group, and 
$\tau_{\le d_x -1}$ denotes the truncation to degrees $\le d_x -1$.
\end{theorem}

(In fact, the Poincar\'e polynomial $P_{G_x/T_x}(t)$ divides
$P_{G/T}(t)$, and the quotient has non-negative coefficients. See
\cite{BP} p.~321.)

To turn these recursive relations into an algorithm for computing the
intersection cohomology Betti numbers, we need a combinatorial
description of all isotropy groups, slices, and links. But such a 
description is unknown in general; in fact, a classification of
spherical homogeneous spaces is only known for $G$ of type $A$, see
\cite{L}.

\medskip

Here we obtain such a combinatorial description for the
subclass of {\it reductive varieties} introduced in \cite{AB1},
\cite{AB2}. It contains all normal, $G\times G$-equivariant embeddings
of the group $G$, and all their invariant subvarieties. Further, both
classes of reductive varieties and of group embeddings are stable
under taking slices and links. Our main tool is the classification of
reductive varieties in terms of certain toric varieties with additional
symmetries, established in [loc.cit.]. The resulting algorithm
specializes to the one in \cite{St1}, \cite{DeLo} and \cite{Fi} for
toric varieties. 

\medskip

The latter algorithm defines remarkable numerical invariants of
rational polytopes, which in fact make sense for non-rational polytopes
as well; see \cite{St1}, \cite{St2}. This has been the starting point for
several recent investigations, constructing a combinatorial
intersection theory for non-rational polytopes, see \cite{BBFK1},
\cite{BBFK2}, \cite{BBFK3}, \cite{BL1}, \cite{BL2}, \cite{Ka}. It
would be interesting to generalize these constructions to the setting
of reductive varieties.

\medskip

The outline of the paper is as follows. Section 2 introduces notation
and basic definitions concerning varieties with algebraic group
actions, and (equivariant) intersection cohomology. In Section 3, we
obtain a variant of a result of \cite{BJ}, that applies to all
spherical varieties, and we study orbits, slices and links in scs
varieties. Theorem \ref{main} is proved in Section 4, and the
combinatorics of reductive varieties are developed in Section 5. The
final Section 6 is devoted to a few examples, including the case of
toric varieties.

\section{Notation and conventions}

\subsection{}\label{conv.1}

We recall the terminology and conventions from our earlier
paper \cite{BJ}. Throughout this paper, we consider complex
algebraic varieties, that is, separated reduced schemes of finite
type over $\bC$. Observe that varieties need not be irreducible;
however, we assume them to be {\it equidimensional}.  

We denote by $G$ a complex linear algebraic group, and by $G^0$ its
connected component containing the identity. A variety $X$ provided
with an algebraic action of $G$ is called a 
$G$-{\it variety}. If, in addition, $X$ admits an equivariant
embedding into the projectivization of a $G$-module, we say that
$X$ is $G$-{\it quasiprojective}. We only consider $G$-varieties
where each orbit admits an open $G$-quasiprojective neighborhood. This 
assumption holds e.g. for normal varieties, see \cite{Su} (for
connected $G$) and \cite{J} (for arbitrary $G$).

\smallskip

\subsection{}\label{slice}

Consider a $G$-variety $X$ and a point $x\in X$; let $Gx$ be its
$G$-orbit and $G_x$ its isotropy group. A {\it slice} to $Gx$ at $x$
is a locally closed subvariety $\cS_x$ of $X$ containing $x$ and
satisfying the following conditions:

\smallskip

(i) $\cS_x$ is invariant under a maximal torus $T_x$ of $G_x$.

\smallskip

(ii) The map $G\times \cS_x \to X, ~(g, s)\mapsto gs$ is 
{\it smooth} at the point $(e,x)$, and the dimension of $\cS_x$ is the
codimension of $Gx$ in $X$.

Such a slice $\cS_x$ always exists, and may be chosen invariant under
a maximal reductive subgroup of $G_x$. Moreover, by shrinking $\cS_x$
if necessary, we may assume that the map $G\times \cS_x \to X$ is 
smooth everywhere, and that $\cS_x$ is affine. 

\subsection{}\label{att.fixed.pt}

Let $T$ denote a torus acting on a variety $X$ with a
fixed point $x$. We say that $x$ is {\it attractive} if there exists a
one-parameter subgroup $\lambda:\bG_m \to T$ such that, for all
$y$ in a Zariski neighborhood of $x$, we have 
$\lim_{t\to 0}\lambda(t)y=x$. Equivalently, all weights of
$T$ acting on the Zariski tangent space at $x$ are contained in an
open half-space.
 
In the situation of (\ref{slice}), we say that $\cS_x$ is an 
{\it attractive slice}, if $x$ is an attractive fixed point for the
action of $T_x$ on $\cS_x$. (See \cite{BJ} (A.1) for further details
on attractive fixed points.) In this case, the geometric quotient
$$
\bP(\cS_x)=(\cS_x - x)/\bG_m
$$ 
exists and we call it the {\it link} at $x$. This is a projective
variety, since $\cS_x$ is assumed to be affine.

\subsection{}\label{coh.1}

Let $X$ be a variety. We denote by $H^*(X)$ the cohomology
ring of $X$ with rational coefficients. $IH^*(X)$ denotes
the corresponding intersection cohomology, for the middle perversity
and rational coefficients. $H^*_G(X)$ ($IH^*_G(X)$) denotes the
corresponding equivariant cohomology ring of $X$ with rational
coefficients (the equivariant intersection cohomology of $X$ with
rational coefficients, \res).

\smallskip

For the sake of completeness we briefly recall the definition of the
intersection cohomology complex (for the middle perversity).
Let $X$ be of dimension $d$ and be provided with a filtration
$U_0 {\overset {j_0} \to} U_1 {\overset {j_1} \to }\cdots  
{\overset {j_{d-2}} \to } U_{d-1} {\overset {j_{d-1}} \to} U_d = X$ 
with each $U_i$ open and each $U_{i}- U_{i-1}$ smooth. Then the
intersection cohomology complex is the complex of locally constant 
sheaves 
$$
IC(X) = \tau_{\le d-1} Rj_{d-1*} \cdots \tau_{\le 0} 
Rj_{0*}({\underline \bQ})
$$
on $X$, where ${\underline \bQ}$ denotes the constant sheaf on $U_0$. 

Assume, in addition, that $X$ is a provided with a $G$-action and that
each $U_i$ is $G$-invariant; let $EG\to BG$ be a universal principal
$G$-bundle. Then 
$EG{\underset G \times} U_i {\overset {j^G_i} \longrightarrow} 
EG{\underset G \times} U_{i+1}$  
provides a filtration of $EG{\underset G \times }X$. Now the
equivariant intersection cohomology complex is 
$$
IC^G(X) = \tau_{\le d-1} Rj_{d-1*}^G \cdots \tau_{\le 0}
Rj_{0*}^G({\underline \bQ}),
$$
where ${\underline \bQ}$ now denotes the constant sheaf on
$EG{\underset G \times} U_0$. The equivariant intersection cohomology
$IH^*_G(X)$ is defined to be the hypercohomology
$\bH^*(EG {\underset G \times} X ; IC^G(X))$. 
These are discussed in more detail in \cite{BJ} Section 1. Both
$H^*_G(X)$ and $IH^*_G(X)$ are graded modules over $H^*(BG)$, the
equivariant cohomology ring of the point.

\smallskip

For any integer $n$, we denote by $\cH^n(IC(X))$ the $n$-th cohomology
sheaf of the intersection cohomology complex on $X$. The stalk
of the sheaf $\cH^n(IC (X))$ at a point $x$ is denoted
$\cH^n(IC(X))_x$, while the local intersection cohomology with support
in $x$ is denoted $IH^n_x(X)$. They are related as follows: 
$IH^n_x(X)$ is the dual space of $\cH^n(IC (X))_x[2d]$, where $d$
denotes the (complex) dimension of $X$.

\section{The key methods}

\subsection{}

We begin by recalling one of the main results of \cite{BJ} (Theorem 2).

\begin{theorem}\label{bj1}
Let $X$ be a $G$-variety containing only finitely many orbits,
each of which admits an attractive slice. Then the following hold.

\smallskip

(i) The $H^*(BG)$-module $IH^*_G(X)$ admits a filtration with
subquotients $IH^*_{\cO,G}(X)$, where $\cO$ runs through the $G$-orbits
in $X$, and $IH^*_{\cO,G}(X)$ denotes the equivariant intersection
cohomology with supports in $\cO$.

\smallskip

(ii) For $\cO=Gx$, the group of components $G_x/G_x^0$ acts on
$H^*(BG_x^0)$ and on $IH^*_x(X)$, and one obtains the isomorphism:
\be\begin{equation}\label{key.1}
IH^*_{\cO,G}(X) \cong IH^{*+2 \dim \cO}_{x, G_x}(X)
\cong (H^{*+2\dim \cO }(BG_x^0) \otimes IH^{*}_x(X))^{G_x/G_x^0}.
\end{equation}\ee
\end{theorem}

One may interpret statement (i) as saying that the stratification 
by orbits is {\it perfect for equivariant intersection cohomology},
under the hypotheses of Theorem \ref{bj1}. However, these hypotheses
are generally not satisfied by spherical varieties. For example, the
linear space $\bC^n$ is spherical under the natural action of the
group $G={\rm SL}_n({\bC})$, and the fixed point $0$ admits no
attractive neighborhood. Likewise, the space of all $n\times n$
matrices of rank at most $1$ is a reductive variety for the same group
(acting by left and right multiplication), but again, the fixed point
$0$ admits no attractive neighborhood. 

For this reason, we will obtain a variant of Theorem \ref{bj1}, where
the existence of attractive slices is replaced by the vanishing of 
local intersection cohomology groups in all odd degrees. The latter 
assumption holds for all spherical varieties by \cite{BJ} Theorem 4;
another proof of that theorem will be given in Remark \ref{alter}. 

\begin{theorem}\label{bj2}
Let $X$ denote a $G$-variety containing only finitely many
orbits. Assume, in addition, that $IH^n_x(X)=0$ for all $x\in X$ and
all odd $n$. Then the conclusions of Theorem \ref{bj1} hold.
\end{theorem}

\begin{proof} 
We first prove (ii). Let $i_\cO:\cO \to X$ denote the inclusion. Then
we obtain 
$$
IH^*_{\cO,G}(X) \cong \bH^*(\cO ; Ri_\cO^! IC^G(X))
\cong \bH^*(BG_x ; Ri_\cO^! IC^G(X)),
$$ 
where the last isomorphism follows from \cite{BJ} (1.6.1). Denoting by
$i_x:x\to X$ the inclusion, we also have 
$$
Ri_\cO^!IC^G (X) \cong Ri_x^!IC^{G_x} (X)[2 \dim (\cO )]
$$
by \cite{BJ} (1.6.2). This yields an isomorphism
$$
IH^*_{\cO,G}(X) \cong 
\bH^{*+2\dim \cO}(BG_x ; Ri_x^!IC^{G_x} (X)).
$$

On the other hand, it follows from \cite{BJ} Theorem 1 that the group
$G_x/G_x^0$ acts on $\bH^*(BG_x^0,Ri_x^!IC^{G_x^0} (X))$, and that
$$
\bH^*(BG_x ; Ri_x^!IC^{G_x} (X)) \cong
\bH^*(BG_x^0 ; Ri_x^!IC^{G_x^0} (X))^{G_x/G_x^0}.
$$
Since $BG_x^0$ is simply connected, the cohomology sheaves of the
complex $Ri_x^!IC^{G_x^0} (X)$ are constant, with
stalks $IH^n_x(X)$. By assumption, these stalks vanish in all odd
degrees; it follows that $Ri_x^!IC^{G_x^0} (X)$ is isomorphic in
$D_b(BG_x^0)$ to its cohomology, $\bigoplus_n IH^n_x(X)[-n]$. This
yields an isomorphism
\be\begin{equation}\label{is}
\bH^*(BG_x^0 ; Ri_x^!IC^{G_x^0} (X)) \cong 
H^*(BG_x^0)\otimes IH^*_x(X).
\end{equation}\ee
Thus, $G_x/G_x^0$ acts on the right-hand side. By \cite{BJ} Lemma 3.6,
it follows that $G_x/G_x^0$ acts on $IH^*_x(X)$ so that the
isomorphism (\ref{is}) is equivariant. This completes the proof of
(ii).

Since $H^*(BG_x^0)$ vanishes in all odd degrees as well, it also
follows that $IH^*_{G,\cO}(X)$ vanishes in all odd degrees. Now choose
$\cO$ closed in $X$, then the long exact sequence
$$
\cdots \to IH^n_{\cO,G}(X) \to IH^n_G(X) \to IH^n_G(X - \cO) \to \cdots
$$
breaks up into short exact sequences. This implies (i) by a
straightforward induction. 
\end{proof}

\subsection{}

Next we review the local structure of spherical varieties. 

Let $G$ be a connected reductive group and $B$ a Borel subgroup with
unipotent radical $U$. We denote by $X$ a spherical $G$-variety, and
by $\rk X$ its {\it rank}, that is, the minimal codimension of a
$U$-orbit; then $\rk X$ is the codimension of $U\xi$, for any point
$\xi$ of the open $B$-orbit.

Choose a $G$-orbit $\cO \subseteq X$; then $\cO$ is spherical as well,
hence we may choose $x\in \cO$ such that $Bx$ is open in $Gx=\cO$. Now
let 
\be\begin{equation}\label{def}
X_0 = \{ \xi \in X ~\vert~ x\in\overline{B\xi}\} = 
\{ \xi \in X~\vert~ \cO\subseteq\overline{B\xi}\}.
\end{equation}\ee
Then $X_0$ is an open affine $B$-invariant subset of $X$, intersecting
$\cO$ along $Bx$. Let $P$ be the normalizer of $X_0$ in $G$.
This is a parabolic subgroup of $G$ containing $B$; let $R_u(P)$ be
its unipotent radical. Now there exists a Levi subgroup $L$ of $P$
and a closed subvariety $\Sigma$ of $X_0$ such that:

\begin{itemize}
\item $\Sigma$ is $L$-invariant and contains $x_0$, and 

\item the map $R_u(P)\times \Sigma \to X_0,~(g,\xi)\mapsto g\xi$
is an isomorphism. 
\end{itemize}

Thus, $\Sigma$ is an affine spherical $L$-variety, of rank equal to
that of $X$. Moreover, $\Sigma\cap Gx$ equals $Lx$; this is the unique
closed $L$-orbit in $\Sigma$. Finally, the isotropy group $P_x$ equals
$L_x$, and contains the derived subgroup $[L,L]$. (See e.g. \cite{K}.)

It follows that $P$ is the normalizer of $Bx$ in $G$. Further, 
$L_x = [L,L] C_x$ where $C$ is the connected center of $L$, and $B_x =
(B\cap [L,L]) C_x$. As a consequence, $L_x^0$ is a connected reductive
group, with Borel subgroup 
$B\cap L_x^0 = (B\cap [L,L]) C_x^0 = B_x^0$, 
of unipotent radical $U\cap L_x^0 = U\cap [L,L] = U_x$.

By a corollary of Luna's slice theorem, there exists a closed
$L_x$-invariant subvariety $\cS_x$ of $\Sigma$, containing $x$, such
that the canonical map
$$
L\times^{L_x} \cS_x \to \Sigma
$$
is an isomorphism. As a consequence, $\cS_x$ is a slice to $Bx$ at $x$,
for the $B$-action on $X$.

\begin{lemma}\label{ranks}
$\cS_x$ is an affine $L_x^0$-spherical variety, and 
$\rk \cS_x = \rk X - \rk Gx$.
\end{lemma}

\begin{proof} 
Since $Lx$ is the unique closed $L$-orbit in $\Sigma$, the point $x$
is the unique closed $L_x$-orbit in $\cS_x$. In particular, $\cS_x$ is
connected.

We claim that $\cS_x$ is normal. To see this, consider the normalization
$\nu:\tilde \cS_x \to \cS_x$. Then the $L_x$-action on $\cS_x$ lifts to an
action on $\tilde \cS_x$ so that $\nu$ is equivariant. Thus, $\nu$ extends
to a morphism
$$
L\times^{L_x}\nu: L\times^{L_x}\tilde \cS_x \to 
L\times^{L_x} \cS_x = \Sigma.
$$
Moreover, the morphism $L\times^{L_x}\nu$ is finite and birational,
since $\nu$ is. But $\Sigma$ is normal, so that $L\times^{L_x}\nu$ is
an isomorphism; thus, the same holds for $\nu$.

Since $\cS_x$ is connected and normal, it is irreducible. And since
$\Sigma=L\times^{L_x} \cS_x$ contains a dense orbit of $B\cap L$, it
follows that $\cS_x$ contains a dense orbit of $B\cap L_x$, and hence of
its subgroup of finite index $B\cap L_x^0$. Thus, $\cS_x$ is a spherical
$L_x^0$-variety; the assertion on ranks follows from the equalities
$\rk X = \rk \Sigma$ and $\rk Gx = \dim Lx = \dim L/L_x$.
\end{proof}

Next we obtain a slight refinement of a result of Knop (\cite{K}
Corollary 7.9 and Remark, p. 326).

\begin{lemma}\label{attr}
There exists an attractive $\bG_m$-action on $\cS_x$ that fixes $x$ and
commutes with the $L_x$-action. 

As a consequence, $L_x$ acts on the link 
$\bP(\cS_x)=(\cS_x - x)/\bG_m$, which is a spherical $L_x^0$-variety
of rank equal to $\rk \cS_x -1$.
\end{lemma}

\begin{proof} 
We use the notation of \cite{K} Section 7. Notice that the sources of
the spherical $L$-variety $\Sigma$ are precisely the closed
$L$-invariant subvarieties; in particular, the closed orbit $Lx$ is a
source. Thus, the closure $\bar{A}_{\Sigma}$ of a generic twisted flat
meets $Lx$. The normalization of $\bar{A}_{\Sigma}$ is an affine
embedding of a finite quotient of the torus $A_{\Sigma}$; let
${\mathcal C}$ be the corresponding cone, then ${\mathcal C}$ is
invariant under the little Weyl group $W_{\Sigma}$. By the argument of
\cite{K} 7.9, there exists a $W_{\Sigma}$-invariant one-parameter
subgroup $v_0$ in the relative interior of ${\mathcal C}$; then both
$v_0$ and $-v_0$ identify to $L$-invariant valuations of the function
field $\bC(\Sigma)$. Thus, $v_0$ yields a $\bG_m$-action on $\Sigma$
commuting with the $L$-action, that is, an $L$-invariant 
grading of the algebra of regular functions $\bC[\Sigma]$. For this
$\bG_m$-action, $\lim_{t\to 0} t \xi$ exists and belongs to
$Lx$ for generic $\xi\in \Sigma$ (since this holds for all
$\xi$ in a generic flat, by definition of $v_0$). It follows that the
corresponding grading of $\bC[\Sigma]$ is non-negative, and that 
$\bC[\Sigma]_0=\bC[Lx]$. Thus, the grading induces a 
positive grading of $\bC[\cS_x]$, the quotient of $\bC[\Sigma]$ 
by the ideal generated by the maximal ideal of $x$ in 
$\bC[Lx]$; this positive grading is clearly $L_x$-invariant.
This proves the first assertion.

>From that assertion and Lemma \ref{ranks}, it follows that
$\bP(\cS_x)$ is a spherical $L_x^0$-variety. To determine its rank,
choose $\xi\in\cS_x$ such that $B_x^0\xi$ is open in $\cS_x$, and let
$[\xi]$ be its image in $\bP(\cS_x)$. Then the isotropy group
$B_{[\xi]}$ acts on the orbit $\bG_m \xi\cong \bC^*$ via a character
with kernel $B_{\xi}$. This implies $U_{[\xi]}=U_\xi$, and hence
$$
\rk \bP(\cS_x) = \dim \bP(\cS_x) - \dim U_x [\xi]
= \dim \cS_x - 1  - \dim U_x \xi = \rk \cS_x -1.
$$ 
\end{proof}

We will also need the following preliminary result.

\begin{lemma}\label{isotrop.1}
Let $\xi$ be a point of the open $B\cap L_x^0$-orbit in $\cS_x$. Then
the orbit $B\xi$ is open in $X$, the isotropy group $B_\xi$ is
contained in $B_x$, and the quotient $B_x/B_\xi$ is irreducible.
\end{lemma}

\begin{proof} By the structure of $X_0$, we have that: $B\xi$ is
open in $X$, and $B_\xi=B\cap L_\xi$. Moreover, since there is a
$L$-equivariant map $\Sigma \to L/L_x$, and since 
$L/L_x = (B\cap L)/(B\cap L_x)$, it follows that 
$B_\xi \subseteq B\cap L_x=B_x$. Note also that the homogeneous space
$B_x/B_\xi$ is the open $B_x$-orbit in $\cS_x$. But $\cS_x$ is
irreducible, so that $B_x/B_\xi$ is irreducible as well.  
\end{proof}

\subsection{}

We may now establish the properties of scs varieties presented in the
Introduction. 

\begin{lemma}\label{isotrop.2}
Let $X$ be a scs $G$-variety. Then all $G$-orbits, slices, and links
in $X$ are scs, and the $G$-isotropy groups of all points are connected.
\end{lemma}

\begin{proof}
Let $\xi\in X$ such that $B\xi$ is open in $X$. Then the product
$BG_{\xi}$ is open in $G$, so that $G_\xi B/B$ is open in
$G/B$, and hence is irreducible. But $G_\xi B/B\cong G_\xi/B_\xi$, and
$B_\xi$ is connected by assumption. Thus, $G_\xi$ is connected as well.

Next consider a $G$-orbit $\cO$ in $X$. Choosing $x$, $\xi$ as in
Lemma \ref{isotrop.1}, the irreducibility of $B_x/B_\xi$ implies that
$B_x$ is connected. Hence $L_x = [L,L]B_x$ is connected as well. 
Thus, the orbit $Gx$ and the slice $\cS_x$ are scs.
It remains to show that the link $\bP(\cS_x)$ is scs. Let $\xi$ as
above ; then, as noted in the proof of Lemma \ref{attr}, 
the isotropy group $B_{[\xi]}$ of the corresponding point of
$\bP(\cS_x)$ is the kernel of a character of $B_{\xi}$. This character
is surjective, since 
$\dim B_{[\xi]} = \dim B - \dim B[\xi] = \dim B -\dim B\xi +1 
= \dim B_{\xi} + 1$. But $B_{\xi}$ is connected, so that $B_{[\xi]}$
is connected. 
\end{proof}

\begin{lemma}
Let $X$ be a scs $G$-variety. Then the $B$-isotropy groups of all
points are connected.
\end{lemma}

\begin{proof}
Let $\xi\in X$. By Lemma \ref{isotrop.2}, we may assume that
$X=G\xi$. We argue by induction on the codimension $c$ of $B\xi$ in
$X$. 

If $c=0$, then $B\xi$ is open in $X$, so that $B_\xi$ is
connected by assumption. If $c\ge 1$, then we may find a minimal
parabolic subgroup $P\supset B$ such that $\overline{B\xi}$ is not
$P$-invariant. Then $B\xi$ is contained in $P\xi$ as a closed $B$-orbit
of codimension $1$, and $P\xi$ contains an open $B$-orbit, say
$B\eta$, of dimension $\dim B\xi + 1$. Thus, $B_\eta$ is connected by
the induction assumption. Further, $P_\eta/B_\eta \cong P_\eta B/B$,
and $P_\eta B$ is open in $P$ (since $B\eta$ is open in $P\eta=P\xi)$,
so that $P_\eta/B_\eta$ is irreducible. Thus, $P_\eta$ is connected as
well.

On the other hand,  the natural map $P\times^B B\xi \to P\xi$ is
finite, since $P$ moves $\overline{B\xi}$. In other words, $B_\xi$ has
finite index in $P_\xi$. But $P_\xi$ is conjugate to $P_\eta$, so that 
$B_\xi$ is connected.
\end{proof}

\section{Proof of the main theorem}

\subsection{}

We begin by introducing various Poincar\'e series. These are formal
power series in a variable $t$, with integer coefficients.  

If $G$ is a linear algebraic group, we put
\be\begin{equation}
P^G(t) = \sum_n \dim H^n(BG) \; t^n,
\end{equation}\ee
the Poincar\'e series of $BG$. For example, if $G$ is a torus
of dimension $r$, then 
\be\begin{equation}
P^G(t) = \frac{1}{(1-t^2)^r}.
\end{equation}\ee
More generally, if $G$ is connected with maximal torus $T$ of
dimension $r$, then we have a fibration $BT\to BG$ with fiber $G/T$
homotopic to the flag manifold of $G$. Hence the cohomology of $G/T$
vanishes in all odd degrees; this implies  
\be \begin{equation}\label{PG}
P^G(t) = \frac{1}{(1-t^2)^r P_{G/T}(t)},
\end{equation}\ee
where $P_{G/T}(t)$ is the Poincar\'e polynomial of $G/T$.

If $X$ is a $G$-variety, we put

\be \begin{equation}
IP^G_X(t) = \sum_n \dim IH^n_G(X) \; t^n.
\end{equation} \ee
In particular, if $G$ is the trivial group, then $IP_X(t)$ is the
intersection cohomology Poincar\'e polynomial of $X$. If $X$ is 
projective of dimension $d$, then we have 
$IP_X(t)=t^{2d} IP_X(\frac{1}{t}),$
by Poincar\'e duality. 

If $G$ is connected and $X$ is a projective $G$-variety, then we have

\be\begin{equation}\label{IPG}
IP^G_X(t) = P^G(t) IP_X(t),
\end{equation}\ee
as follows from the degeneration of the spectral sequence in
equivariant intersection cohomology (see e.g. \cite{BJ} (1.5.2)).

Returning to an arbitrary $G$-variety $X$, we put for any $x\in X$:
\be\begin{equation}
IP^{G_x}_{x,X}(t) = \sum_n \dim IH^n_{x,G_x}(X) \; t^n.
\end{equation}\ee

In particular, $IP_{x,X}(t) =\sum_n \dim IH^n_x(X) \; t^n$ is the
Poincar\'e polynomial for local intersection cohomology with support
in $x$.

On the other hand, one also has the Poincar\'e series for the stalks
of the equivariant intersection cohomology sheaves:
\be\begin{equation}
IP_{X,x}^{G_x}(t) = \sum_n \dim \cH^n(i_x^*IC^{G_x}(X)) \; t^n,
\end{equation}\ee
where $i_x: x \to X$ denotes the inclusion. 

Note that the Poincar\'e polynomial 
$IP_{X,x}(t) = \sum_n \dim \cH^n(i_x^*IC(X)) \; t^n$
satisfies by Verdier duality :
\be\begin{equation}\label{IP}
IP_{X,x}(t) = t^{2d} IP_{x,X}(\frac{1}{t}).
\end{equation}\ee

\subsection{}

We may now formulate a direct consequence of Theorems \ref{bj1} and
\ref{bj2} for these Poincar\'e series.

\begin{proposition}\label{rec}

\smallskip 

\noindent
(i) Let $X$ be a $G$-variety satisfying the assumptions of Theorem
\ref{bj1} or \ref{bj2}. Then 
$$
IP^G_X(t) = \sum_x t^{-2\dim Gx} \; IP^{G_x}_{x,X}(t)
$$
(sum over representatives of $G$-orbits in $X$).

\smallskip

\noindent
(ii) If, in addition, $G_x$ is connected, then 
$$
IP^{G_x}_{x,X} (t) = P^{G_x}(t) IP_{x,X}(t).
$$
\end{proposition}

Next let $X$ be a projective scs variety. Then Proposition
\ref{rec} and (\ref{PG}), (\ref{IPG}), (\ref{IP}) imply
$$
IP_X(t) = \sum_x (1-t^2)^{r-r_x}
\frac{P_{G/T}(t)}{P_{G_x/T_x}(t)}t^{2d_x} IP_{X,x}(\frac{1}{t}).
$$
Replacing $t$ with $\frac{1}{t}$ and using Poincar\'e duality for $X$,
$G/T$ and $G_x/T_x$ yields the equality (\ref{A.1}) in Theorem
\ref{main}. Now the equality (\ref{A.2}) is a consequence of the
following statement.

\begin{proposition}\label{trunc}
Consider a $G$-variety $X$, a point $x\in X$, and a slice $\cS_x$ to
$Gx$ at $x$. Then 
\be\begin{equation}\label{sl}
IP_{x,X}(t) = t^{2(d-d_x)} IP_{x,\cS_x}(t)
= t^{2d} IP_{\cS_x,x}(\frac{1}{t}),
\end{equation}\ee
where $d=\dim X $ and $d_x=\codim Gx = \dim \cS_x$. 

If, in addition, $\cS_x$ is attractive with link $\bP(\cS_x)$, then 
\be\begin{equation}\label{tr}
IP_{\cS_x,x}(t) = \tau_{\le d_x-1}((1-t^2)IP_{\bP(\cS_x)}(t))
\end{equation}\ee
if $d_x \ge 2$, and $IP_{\cS_x,x}(t)=1$ otherwise.
\end{proposition}

\begin{proof}
One obtains a quasi-isomorphism  
$$
Ri_{x,X}^!IC(X) \cong Ri_{x,\cS_x}^!IC(\cS_x)[-2(d-d_x)],
$$ 
where $i_{x,X}: x \to X$ and $i_{x,\cS_x}:x \to \cS_x$ denote the
inclusions. This yields the first equality in (\ref{sl}); the second one
follows from (\ref{IP}).

For (\ref{tr}), observe that 
$\cH^n(IC(\cS_x))_x \cong IH^n(\cS_x) \cong IH^n(\cS_x-x)$ for 
all $n \le d_x-1$, while $\cH^n(IC(\cS_x))_x =0$ for all $n>d_x-1$ 
(see, for example, the definition of the intersection cohomology
complex in \ref{coh.1}). Thus, 
$IP_{\cS_x,x}(t) = \tau_{\le d_x -1}(IP_{\cS_x - x}(t))$. Now consider
the Wang exact sequence
$$\displaylines{
\cdots \to IH^{n-2}(\bP(\cS_x)) \to IH^n(\bP(\cS_x)) \to IH^n(\cS_x-x)
\to 
\hfill\cr\hfill
\to IH^{n-1}(\bP(\cS_x)) \to IH^{n+1}(\bP(\cS_x)) \to \cdots,
\cr}$$
where the first and last maps are multiplication by the class of a
hyperplane in $\bP(\cS_x)$ (see e.g. \cite{BJ} 3.5). By the Hard
Lefschetz theorem in intersection cohomology, these maps are injective 
for all $n \le \dim \bP(S_x) = \dim \cS_x-1 = d_x -1$. Therefore, the
Wang exact sequence breaks up into short exact sequences for all 
$n \le d_x-1$. This yields 
$\tau_{\le d_x -1}(IP(\cS_x - x)) = 
\tau_{\le d_x -1}((1-t^2)IP_{\bP(\cS_x)}(t))$.
\end{proof}

\begin{remark}\label{alter}
The results of Sections 3 and 4 yield another proof of the vanishing
of $IH^n(X)$ and $IH^n_x(X)$ for all odd $n$ and spherical $X$
(\cite{BJ} Theorem 4). Indeed, let us argue by induction on
$\dim X = d$. By Lemma \ref{attr}, all links of $X$ are spherical, of
dimension $<d$. Thus, by Proposition \ref{trunc} and the induction
assumption, $IH^n_x(X)$ for all odd $n$ and $x\in X$. Now Theorem
\ref{bj2} implies the vanishing of $IH^n(X)$ for all odd $n$.
\end{remark}

\begin{remark}\label{virtu}
For any closed connected subgroup $H\subseteq G$, put 
\be\begin{equation}
P_{G/H}(t) = (-1)^{r-s} \frac{P^H(t)}{P^G(t)},
\end{equation}\ee
where $s$ is the rank of $H$ (and $r$ is the rank of $G$). 
Then $P_{G/H}(t)$ is the {\it virtual Poincar\'e polynomial of $G/H$},
as follows from \cite{DiLe} Theorem 6.1 (ii) applied to the fibration
$BH\to BG$ with fiber $G/H$. In particular, $P_{G/H}(t)$ is a
polynomial with rational coefficients, of degree $\dim G/H$, that
only depends on the structure of complex algebraic variety on $G/H$. 
(See \cite{BP} for more on virtual Poincar\'e polynomials of
homogeneous spaces.) 

With this notation, (\ref{PG}) yields
\be\begin{equation}
P_{Gx}(t) = P_{G/G_x}(t) = (t^2 - 1)^{r-r_x} 
\frac{P_{G/T}(t)}{P_{G_x/T_x}(t)},
\end{equation}\ee
so that (\ref{A.1}) translates into 
\be\begin{equation}
IP_X(t)= \sum_x P_{Gx}(t)\, IP_{X,x}(t).
\end{equation}\ee
We may regard the terms in the right-hand side as virtual Poincar\'e
polynomials for intersection cohomology with supports in orbits.
\end{remark}

\section{Reductive varieties}

\subsection{}

We first give an overview of {\it affine} reductive varieties, after
\cite{AB1}. For this, we need to introduce notation concerning
reductive groups. 

Let $G$ be a connected reductive group and let $B$, $B^-$ be opposite
Borel subgroups, i.e., $T=B\cap B^-$ is a maximal torus of $G$; let
$U=R_u(B)$, $U^-=R_u(B^-)$ be the unipotent radicals of $B$, $B^-$.
The character group of $T$ is denoted by $\Lambda$, and called the
weight lattice; we put $\Lambda_{\bR}=\Lambda\otimes_{\bZ}{\bR}$.
Let $W$ be the Weyl group of $(G,T)$; it acts on $\Lambda$ and on
$\Lambda_{\bR}$. The root system of $(G,T)$ is denoted by $\Phi$, 
with the subsets $\Phi^+$ of positive roots (the roots of $(B,T)$),
and $\Pi$ of simple roots. 

For a subset $I\subseteq \Pi$, we denote by $\Phi_I$ the corresponding
subsystem of $\Phi$, and by $W_I$ (resp. $P_I\supseteq B$ and
$P^-_I\supseteq B^-$) the corresponding parabolic subgroup of $W$
(resp. opposite parabolic subgroups of $G$). We put 
$L_I=P_I\cap P_I^-$; this is a Levi subgroup of $P_I$ and 
$P_I^-$, with root system $\Phi_I$ and Weyl group $W_I$. Let $\ell$ be
the length function of $W$, and let $W^I$ be the subset of
representatives of minimal length of $W/W_I$. Then the Poincar\'e
polynomial of $G/P_I$ equals $\sum_{w\in W^I} t^{2\ell(w)}$.   
 
Now consider the connected reductive group $G\times G$, with Borel
subgroup $B^-\times B$ and maximal torus $T\times T$. Let $\diag T$ be
the diagonal in $T\times T$. We say that an affine $G\times G$-variety
$X$ is {\it reductive} (for $G$) if it satisfies the following
conditions.

\smallskip

\noindent
(i) $X$ is normal.

\smallskip

\noindent
(ii) There exists $x\in X$, fixed by $\diag T$, such that the orbit
$(B^-\times B)x$ is dense in $X$.

\smallskip

\noindent
(iii) The isotropy group $(G\times G)_x$ is connected.

\smallskip

Further, we may replace (iii) with the assumption of connectedness of
$(B^-\times B)_x$, or of $(T\times T)_x$. Thus, affine reductive
varieties are scs $G\times G$-varieties. Note also that the set of all
$x\in X$ satisfying (ii) is a unique $T\times T$-orbit. Any such point
is called a {\it base point}. Also note that 
$\dim (T\times T) x = \rk X$. 

By the Bruhat decomposition, the multiplication of $G$ yields an open
immersion $U^-\times T\times U \to G$. Thus, the group $G$, regarded as
a $G\times G$-variety via left and right multiplication, is an affine
reductive variety. More generally, all affine $G\times G$-equivariant
embeddings of $G$, or of quotients of $G$ by connected normal
subgroups, are reductive varieties for $G$. 

Next we summarize results of \cite{AB1} about the classification and
orbit structure of affine reductive varieties. Let $X'$ be the closure
in $X$ of the orbit of base points; this is an affine toric variety for
$T$ (identified with $T\times \{1\} \subset G\times G$). Thus, $X'$ is
uniquely determined by the set of weights of $T$ in its coordinate
ring. Further, this set is the intersection of $\Lambda$ with a
uniquely determined rational convex polyhedral cone $\sigma$ in
$\Lambda_{\bR}$.

The assignement $X\mapsto \sigma$ yields a bijective correspondence
from affine reductive varieties to rational convex polyhedral
cones $\sigma\subseteq \Lambda_{\bR}$ satisfying the following
conditions: 

\smallskip

\noindent
(i) The relative interior $\sigma^0$ meets $\Lambda^+_{\bR}$. 

\smallskip

\noindent
(ii) The distinct $w\sigma^0$ ($w\in W$) are disjoint.

\smallskip

We then say that $\sigma$ is a {\it $W$-admissible cone}, and put
$X=X_\sigma$; then $\rk X = \dim \sigma$. The $G\times G$-orbit
closures in $X_\sigma$ are the affine reductive varieties associated
to $W$-admissible faces of $\sigma$. 

For any $W$-admissible cone $\sigma$, let $N_W(\sigma)$
(resp. $C_W(\sigma)$) be its normalizer (resp. centralizer) in $W$. 
Then we have $N_W(\sigma)=W_I$ and $C_W(\sigma)=W_J$ for subsets 
$J=J(\sigma)\subseteq I=I(\sigma)\subseteq \Pi$. 
Further, all roots in the complement $K = K(\sigma) = I -J$ are
orthogonal to $J$, so that $W_I = W_J\times W_K$ and $L_I = L_J L_K$;
and $J$ (resp. $K$) consists of those simple roots such that the
corresponding wall of $\Lambda^+_{\bR}$ contains $\sigma$ (resp. meets
$\sigma^0$).

Let $\Lambda_\sigma$ be the subgroup of $\Lambda$ spanned by all
elements of $\Lambda\cap\sigma$; this is a saturated sublattice of
$\Lambda$, that is, the quotient $\Lambda/\Lambda_\sigma$ 
is torsion-free. Let $T_\sigma\subseteq T$ be the intersection of the
kernels of all characters in $\Lambda_\sigma$ (or, equivalently, in
$\Lambda\cap\sigma$); this is a central 
subtorus of $L_K$, with character group $\Lambda/\Lambda_\sigma$. 
Put $G_\sigma = [L_J,L_J] T_\sigma$; this is a connected reductive
subgroup of $L_J$ (denoted by $H_J$ in \cite{AB1}), with maximal torus
$T_\sigma$. Thus, $G_\sigma$ is a normal subgroup of $L_I$, and the
quotient $L_I/G_\sigma$ is isomorphic to $L_K/T_\sigma$. Now
\be\begin{equation}\label{isot}
(G\times G)_x = 
(R_u(P_I)\times R_u(P_I^-))(G_\sigma\times G_\sigma) \diag(L_K),
\end{equation}\ee
up to conjugacy by an element of $T\times T$. 

As a consequence, any $G\times G$-orbit in an affine reductive variety
admits a homogeneous fibration over a product of two opposite flag
manifolds, with fiber a connected reductive group. Specifically, we
have a $G\times G$-equivariant morphism 
\be\begin{equation}\label{base}
(G\times G)x \to G/P_I \times G/P_I^-
\end{equation}\ee
with fiber 
\be\begin{equation}\label{fiber}
(L_I\times L_I)/(G_\sigma\times G_\sigma)\diag L_K = 
L_I/G_\sigma = L_K/T_\sigma.
\end{equation}\ee 
Here $P_I\times P_I^-$ acts on this fiber via its quotient 
$L_I\times L_I$ acting on $L_I/G_\sigma$ by left and right
multiplication. Note that $L_K/T_\sigma$ is a connected reductive
group with weight lattice $\Lambda_\sigma$ and root system $\Phi_K$.

It follows readily that 
\be\begin{equation}\label{dim}
\dim X_\sigma = \dim(G\times G)x =
\vert \Phi - \Phi_J \vert + \dim \sigma, 
\end{equation}\ee
and that the virtual Poincar\'e polynomial of $(G\times G)x$ (see
Remark \ref{virtu}) is given by 
\be\begin{equation}\label{pol}
P_{(G\times G)x}(t) = (t^2-1)^{\dim \sigma}
(\sum_{w\in W^I} t^{2\ell(w)})^2 (\sum_{w\in W_K} t^{2\ell(w)}).
\end{equation}\ee

As another consequence of (\ref{isot}), $X$ is an affine
embedding of the quotient of $G$ by a connected normal subgroup if and
only if $I=\Pi$, that is, $\sigma$ is $W$-invariant. In particular,
affine embeddings of $G$ correspond to $W$-invariant cones with
nonempty interior. 

\subsection{}

Following \cite{AB2} Section 2, we define {\it projective} reductive
varieties and we sketch how to deduce their main properties from the
affine case.

Consider a projective irreducible $G\times G$-variety $X$ equipped
with an ample $G\times G$-linearized line bundle $L$. Let 
$R = \bigoplus_{n=0}^{\infty} \Gamma (X,L^{\otimes n})$, 
this is a graded, finitely generated reduced algebra, where 
$G\times G$ acts. This defines an affine variety $\tilde X$ where
$\bG_m\times G\times G$ acts. Further, the action of $\bG_m$ is
attractive, and the corresponding link is nothing but $X$. We say that
the pair $(X,L)$ is a {\it linearized projective $G\times G$-variety}. 

Put $\tilde G = \bG_m\times G$, this is a connected reductive group
with weight lattice $\tilde \Lambda = \bZ \times \Lambda$. We may
regard $\tilde X$ as a $\tilde G\times \tilde G$-variety, where
$\bG_m\times \bG_m$ acts via its morphism 
$(t_1,t_2)\mapsto t_1 t_2^{-1}$ to $\bG_m$.
For any $x\in X$ with representative $\tilde x\in \tilde X$, we
obtain readily an exact sequence of isotropy groups:
\be\begin{equation}\label{exa}
1 \to \bG_m \to (\tilde G\times\tilde G)_{\tilde x} \to
(G\times G)_x \to 1.
\end{equation}\ee

We say that $X$ is reductive for $G$, if $\tilde X$ is reductive 
for $\tilde G$; then $(X,L)$ is called a 
{\it linearized reductive variety}. These may be characterized as
those linearized projective $G\times G$-varieties $(X,L)$ that satisfy
the following conditions:

\smallskip

\noindent
(i) $X$ is normal.

\smallskip

\noindent
(ii) There exists $x\in X$, fixed by $\diag T$, such that the orbit
$(B^-\times B)x$ is dense in $X$, and that $\diag T$ fixes the fiber of
$L$ at $x$.

\smallskip

\noindent
(iii) The isotropy group $(G\times G)_x$ is connected.

\smallskip

Again, (iii) is equivalent to the assumption of connexity of 
$(B^-\times B)_x$, or of $(T\times T)_x$; and the set of all $x\in X$
satisfying (ii) is a unique $T\times T$-orbit: the orbit of base
points, of dimension equal to $\rk X$.

\smallskip
Thus, projective reductive varieties are scs. Again, examples include 
$G\times G$-equivariant embeddings of the quotient of $G$ by a
connected normal subgroup.

Let $\sigma\subseteq \tilde\Lambda_{\bR} = \bR\times \Lambda_{\bR}$
be the cone associated to $\tilde X$, and put 
$\delta=\sigma\cap (1\times \Lambda_{\bR})$. Then $\delta$ is a
lattice polytope in $\Lambda_{\bR}$, and $\sigma$ is the cone over
$\delta$. Since $\sigma$ is $W$-admissible, $\delta$ satisfies the
following conditions:

\smallskip

\noindent
(i) The relative interior $\delta^0$ meets $\Lambda^+_{\bR}$.

\smallskip

\noindent
(ii) The distinct translates $w\delta^0$ ($w\in W$) are disjoint.

\smallskip

A lattice polytope $\delta\subset\Lambda_{\bR}$ satisfying (i) and
(ii) is called a {\it $W$-admissible polytope}. These classify
polarized reductive varieties; we denote by $(X_\delta,L_\delta)$ the
linearized reductive variety with polytope $\delta$, then 
$\dim X_\delta = \rk \delta$. The closure
in $X$ of the orbit of base points, equipped with the restriction of
$L$, is the linearized toric variety with polytope $\delta$. The
$G\times G$-orbit closures in $X_\delta$ are the $X_\phi$, where
$\phi\subseteq \delta$ is a $W$-admissible face.

Since $N_W(\sigma)=N_W(\delta)$ and $C_W(\sigma)=C_W(\delta)$, we
obtain two subsets $J=J(\delta)\subseteq I=I(\delta)\subseteq \Pi$ 
satisfying the properties of the previous subsection.
Now the description (\ref{isot}) of the isotropy group $(G\times G)_x$
carries over to this projective setting, with $\Lambda_\sigma$ being
replaced by the lattice $\Lambda_\delta$ spanned by the 
{\it differences} of any two elements of $\Lambda\cap\delta$. 

As a consequence, the description of orbits as fibered spaces carries
over as well; specifically, the analogues of (\ref{base}),
(\ref{fiber}), (\ref{dim}) and (\ref{pol}) hold with $\sigma$ being
replaced by $\delta$. Further, projective embeddings of a quotient of
$G$ by a connected normal subgroup (resp. of $G$) corrrespond to
$W$-invariant lattice polytopes (resp. with nonempty interior). 

\subsection{}

We obtain a combinatorial description of slices and links in
reductive varieties. 

Consider a $W$-admissible polytope $\delta\subset\Lambda_\bR$, and a
$W$-admissible face $\varphi \subseteq \delta$. These correspond to a
linearized reductive variety $(X_\delta,L_\delta)$ together with a
$G\times G$-orbit $\cO = \cO_\varphi$: the open orbit in 
$X_\varphi\subseteq X_\delta$. We describe the local structure of $X$
along $\cO$, by making explicit the objects introduced in 3.2. 

Let $x$ be a base point of $\cO$, then $(B^-\times B)x$ is open in
$(G\times G)x=\cO$. Further, it follows from (\ref{isot}) that the
normalizer $P$ of $(B^-\times B)x$ in $G\times G$ equals 
$P_J^-\times P_J$, where $J=J(\varphi)$. Since $x$ is fixed by 
$\diag T$, the Levi subgroup $L$ of $P$ equals 
$L_J\times L_J$. Further, by \cite{AB1} Lemma 2.8, the variety
$\Sigma$ is an affine reductive variety for $L_J$ ; one readily checks
that the corresponding $W_J$-admissible cone is generated by the
differences $\lambda - \mu$, where $\lambda\in\delta$ and
$\mu\in\varphi$.

Now by (\ref{isot}) again, we have  $L_x = G_\varphi \times G_\varphi$. 
Note that $G_\varphi$ is a connected reductive subgroup of $G$,
normalized by $T$; further, $T_\varphi$ is a maximal torus of
$G_\varphi$, so that the weight lattice of $G_\varphi$ equals
$\Lambda_\varphi=\Lambda/\Lambda_\varphi$. The set of simple roots of
$G_\varphi$ is $J=J(\varphi)$, with Weyl group $W_J=C_W(\varphi)$; we
denote the latter by $W_\varphi$.

By \cite{AB1} Lemma 4.1, the slice $\cS_x$ is an affine reductive
variety for $G_\varphi$. Denote its $W_\varphi$-admissible cone by
$\sigma = \sigma_\varphi$; this cone is the image in
$\Lambda_\varphi$ of the cone of $\Sigma$. So we may regard $\sigma$ as
the normal cone to $\delta$ along its face $\varphi$. Note the
equality $\rk \cS_x = \dim \delta - \dim \varphi$.

To describe the link $\bP(\cS_x)$, note first that the closed convex
cone $\sigma$ contains no line. Thus, we may find a linear form $f$ on
$\Lambda_{\bR}/\Lambda_{\varphi,\bR}$ that takes positive values at
all non-zero points of $\sigma$. We may assume, in addition, that $f$
takes integer values at all points of $\Lambda/\Lambda_\varphi$, and
is invariant under the normalizer of $\sigma$ in $W_\varphi$. Then by
\cite{AB1} 3.2, 4.1, $f$ yields a positive $G_\varphi\times
G_\varphi$-invariant grading of the algebra of regular functions on
$\cS_x$. In other words, $f$ defines an attractive $\bG_m$-action on
$\cS_x$ that commutes with the action of 
$G_\varphi \times G_\varphi$. Now $\bP(\cS_x)$ is  
the reductive variety for $G_\varphi$ associated with the polytope 
$\sigma\cap (f=n)$, where $n$ is a suitable positive integer. We may
regard this polytope as the link of $\delta$ along its face $\varphi$;
we have $\rk \bP(\cS_x) = \dim \delta - \dim \varphi -1$.

If $X_\delta$ is an embedding of a quotient of $G$ by a connected
normal subgroup, then $\delta$ is $W$-invariant, so that
$\sigma_\varphi$ is invariant under $W_\varphi$. Thus, $\cS_x$
is an embedding of a quotient of $G_\varphi$ by a connected normal
subgroup. So the class of embeddings of connected reductive groups is
stable under taking slices and, likewise, links.

\section{Examples}

\subsection{}

We begin with the case of {\it toric varieties}, where the objects
introduced in Section 5 take a very simple form. 

Let $X$ be a projective toric variety; let $T$ be the corresponding
torus, with character group $\Lambda$. Then $X$ corresponds to a fan
$\Sigma$, consisting of the normal cones to all faces of a lattice
polytope $\delta\in\Lambda_{\bR}$. (This polytope is not uniquely
determined by $X$, but so are the partially ordered set of its faces,
and their directions.)

The $T$-orbits in $X$ correspond to cones of $\Sigma$, that is, to
faces of $\delta$. For any such face $\varphi$, the isotropy group of
the corresponding orbit $\cO = \cO_{\varphi}$ is a subtorus
$T_\varphi$ of $T$, the intersection of all characters in the space
$\Lambda_{\varphi,\bR}$ spanned by the differences $\lambda-\mu$ where
$\lambda,\mu\in\varphi$. 

Moreover, $\cO$ admits an open $T$-invariant neighborhood in
$X$, isomorphic to the product $\cO\times \cS_{\varphi}$,
where $S_{\varphi}$ is an affine toric variety for $T_\varphi$ with a
fixed point $x_{\varphi}$. The cone associated with $\cS_{\varphi}$ is
dual to the cone $\sigma_\varphi$, image in
$\Lambda_{\bR}/\Lambda_{\varphi,\bR}$ of the cone generated by the
differences $\lambda-\mu$ where $\lambda\in\delta$, $\mu\in\varphi$.

It follows that $\cS_{\varphi}$ is an attractive slice at any
point of $\cO$, and that the associated link 
$\bP(\cS_{\varphi})$ is a projective toric variety with polytope the
link of $\varphi$ in $\delta$. The $x_\varphi$, $\varphi$ a face of
$\delta$, form a system of representatives of the $T$-orbits in
$X$. (See e.g. \cite{Fu} 1.4, 2.1.)

Now Theorem \ref{main} yields the equalities

\be\begin{equation}\label{lg}
IP_X(t)= \sum_\varphi (t^2-1)^{\dim \varphi} IP_{X,x_{\varphi}}(t)
\end{equation}\ee
(sum over all faces of $\delta$), and

\be\begin{equation}\label{gl}
IP_{X,x_{\varphi}}(t) = IP_{\cS_{\varphi},x_{\varphi}}(t)
= \tau_{\le \codim\varphi -1} ((1-t^2)IP_{\bP(\cS_{\varphi})}(t)).
\end{equation}\ee

When expressed in terms of cones of $\Sigma$, (\ref{gl}) and
(\ref{lg}) give back the main result of \cite{Fi}; see also
\cite{DeLo}, \cite{St1}, \cite{St2}.

\subsection{}

Next we describe orbits, slices and links in 
{\it reductive varieties of rank $1$.} Indeed, the rank, rather than
the dimension, measures how complicated a spherical variety is; and
spherical (resp. reductive) varieties of rank $0$ are just flag
manifolds $G/P$ (resp. products of two opposite flag manifolds
$G/P_I\times G/P_I^-$). 

By 5.2, a linearized reductive variety $(X,L)$ of rank $1$ corresponds
to a $W$-admissible line segment $\delta = [\lambda,\mu]$, where 
$\lambda,\mu\in\Lambda$. We may assume that $\lambda$ is dominant;
then it is a $W$-admissible face, corresponding to a closed 
$G\times G$-orbit $\cO=\cO_\lambda$. With the notation of 5.2, we have 
$\Lambda_\lambda=0$, so that $T_\lambda=T$. Further,
$I(\lambda)=J(\lambda)$ is the set of simple roots orthogonal to
$\lambda$, so that $G_\lambda=L_{I(\lambda)}$ is the common Levi
subgroup to $P_{I(\lambda)}$ and $P_{I(\lambda)}^-$. Moreover, 
$\cO_\lambda = G/P_{I(\lambda)}\times G/P_{I(\lambda)}^-$. 
The corresponding slice is the affine reductive variety for
$L_{I(\lambda)}$ with cone the ray spanned by $\mu-\lambda$. Setting 
$I(\delta)=I(\lambda)\cap I(\mu)$, the corresponding link is the
product
$$
L_{I(\lambda)}/P_{I(\delta})\cap L_{I(\lambda)}
\times 
L_{I(\lambda)}/P_{I(\delta})^-\cap L_{I(\lambda)}
\cong 
P_{I(\lambda)}/P_{I(\delta)} 
\times 
P_{I(\lambda)}^-/P_{I(\delta)}^-.
$$

If $\mu$ is also dominant, then $X$ contains three 
$G\times G$-orbits: the open orbit, and the closed orbits
corresponding to $\lambda$ and $\mu$. The isotropy group of a base
point of the open orbit $\cO_\delta$ equals 
$(R_u(P_{I(\delta)})\times R_u(P_{I(\delta)}^-))
(G_\delta\times G_\delta)$,
where $G_\delta\subset L_{I(\delta)}$ is the connected kernel of
the character $\mu-\lambda$. Thus, $\cO_\delta$ fibers over 
$G/P_{I(\delta)}\times G/P_{I(\delta)}^-$ with fiber $\bG_m$.

But if $\mu$ is not dominant, then $\lambda$ is the unique
$W$-admissible proper face, so that $X$ contains only two orbits. In
that case, $\lambda$ and $\mu$ are exchanged by a unique simple
reflection $s_\alpha$, where $\alpha\in\Pi$, so that $\mu-\lambda$ is
a non-zero multiple of $\alpha$. Moreover, the isotropy
group of a base point of $\cO_\delta$ equals
$(R_u(P_I)\times R_u(P_I^-)) (G_\delta\times G_\delta)\diag L_\alpha$,
where $I=I(\delta)\cap\{\alpha\}$, and $G_\delta\subset L_{I(\delta)}$
is the connected kernel of the character $\mu-\lambda$ (or of
$\alpha$). Thus, $\cO_\delta$ fibers over 
$G/P_{I(\delta)}\times G/P_{I(\delta)}^-$ with fiber 
$L_\alpha/\ker(\alpha)^0$, isomorphic to $\SL_2$ or $\PGL_2$.

\subsection{}

Finally, as examples of reductive varieties of rank $2$, we consider
{\it embeddings of the group $\GL(2)$.} Let $B$, $B^-$ be
the opposite Borel subgroups of upper, \res lower triangular matrices,
then $T=B\cap B^-$ is the maximal torus of diagonal matrices. Via the
diagonal coefficients, the weight lattice $\Lambda$ identifies to
$\bZ^2$; then the unique simple root is $\alpha=(1,-1)$, the unique
non-trivial element of the Weyl group is the reflection with respect
to the diagonal of $\bR^2$, and the positive Weyl chamber consists of
all points below the diagonal.

Thus, projective embeddings of $\GL(2)$ correspond to lattice
polygons in $\bR^2$, symmetric with respect to the diagonal. Further,
orbit closures in such embeddings correspond to the following four
types of $W$-admissible faces $\varphi$:

\smallskip

\noindent
(i) $\varphi$ is an edge entirely below the diagonal. Then both sets 
$J(\varphi)$, $K(\varphi)$ are empty, so that
$G_\varphi=T_\varphi=T$. The isotropy group of a base point of the
orbit $\cO_\varphi$ is the product of $U\times U^-$ with a torus of
dimension $3$. The corresponding slice is an affine line, and the link
is just a point.

\smallskip

\noindent
(ii) $\varphi$ is an edge, symmetric with respect to the
diagonal. Then $J(\varphi)=\emptyset$, $K(\varphi)=\{\alpha\}$, and  
$G_\varphi=T_\varphi$ is a torus of dimension $1$. The isotropy group
of $\cO_\varphi$ is the product of $\diag \GL(2)$ with
$T_\varphi\times 1$. Again, the slice is a line, and the link is
a point.

\smallskip

\noindent
(iii) $\varphi$ is a vertex below the diagonal. Then
$J(\varphi)=K(\varphi)=\emptyset$, $G_\varphi=T_\varphi=T$, and 
the isotropy group of $\cO_\varphi$ is $B\times B^-$. So
this orbit is isomorphic to $\bP^1\times \bP^1$. The slice is an
affine toric variety of dimension $2$, so that the link is $\bP^1$.

\smallskip

\noindent
(iv) $\varphi$ is a vertex on the diagonal. Then
$J(\varphi)=\{\alpha\}$, $K(\varphi)=\emptyset$, $T_\varphi=T$, and
$G_\varphi =\GL(2)$. The isotropy group of $\cO_\varphi$ is the whole
$\GL(2)\times\GL(2)$, so that $\cO_\varphi$ is a fixed point. The
slice is the affine embedding of $\GL(2)$ associated with the tangent
cone to $\delta$ at its vertex $\varphi$. The corresponding link is a
projective embedding of the quotient of $\GL(2)$ by a non-trivial
central torus, that is, of $\PGL(2)$. It follows that this link is the
projective space $\bP^3$, the projectivization of the space of
$2\times 2$ matrices.


\end{document}